\documentclass[12pt,psfig]{article}
\usepackage{graphicx,epsfig}
\usepackage{mathrsfs}
\usepackage{amssymb}
\def\bbb{\begin{eqnarray*}}

\def\eee{\end{eqnarray*}}

\begin{document}

\baselineskip=17pt

\begin{center}

\vspace{-0.6 in} {\large \bf  On the stability of self-adjointness of linear relations}
\\ [0.3in]

Yan Liu\\ 

Department of Mathematics and Physics, Hohai University,\\
Changzhou Campus 213022, P. R. China

\footnote{The corresponding author.}
\footnote{Email addresses: yanliumaths@126.com}
\end{center}

{\bf Abstract.}
This paper focuses on the stability of self-adjointness of linear relations under perturbations in Hilbert spaces.
It is shown that a self-adjoint relation is still self-adjoint under bounded and relatively bounded perturbations.
The results obtained in the present paper generalize the corresponding results
for linear operators to linear relations,
and some weaken the conditions of the related existing results.
\medskip

\noindent{\bf 2010 AMS Classification}: 47A06, 47A55, 47B25.
\medskip

\medskip
\parindent=10pt

\section{ Introduction }

Perturbation theory are one of the main topics in both pure and applied mathematics.
Since the self-adjoint operators form the most important class of linear operators that appear in applications,
the perturbation of self-adjoint operators and the stability of self-adjointness have received lots of attention.
In particular, Kato first studied the stability of self-adjointness of closed symmetric operators and
showed that the self-adjointness is preserved under relatively bounded perturbations with relative bounds less than 1 [6].
Followed by this work, Devinatz, Zettl, Behncke, Kissin, etc.
extended this work about stability of self-adjointness to results about the stability of the deficiency indices [2, 4, 7, 17].

With further research of operator theory, more and more multi-valued operators and non-densely defined operators have been found.
For example, the operators generated by those linear continuous Hamiltonian systems,
which do not satisfy the definiteness conditions, and general linear discrete Hamiltonian systems
may be multi-valued or not densely defined in their corresponding Hilbert spaces (cf. [8, 10, 14]).
So the classical perturbation theory of linear operators is not available in this case.
Motivated by this need,
von Neumann [9] first introduced linear relations into functional analysis,
and then Arens [1] and many other scholars further studied and developed the fundamental theory of linear relations.
A liner relation is also called a linear subspace (briefly, subspace).

Since the theory of linear relations was established, the related perturbation problems have attracted extensive attention of scholars
and some good results have been obtained [12, 15, 16].
It is well known that the self-adjoint relations are the most important class of linear relations that appear applications.
To the best of our knowledge, there seem a few results about the stability of self-adjointness of linear relation under perturbations [12, 16].
But it has not been specifically and thoroughly studied. In the present paper, we shall concentrate on the stability of self-adjointness of linear relations.
The results given in the present paper not only weaken the conditions of Theorem 4.1 in [12],
but also cover the result obtained in [16, Theorem 5.2].

The rest of the paper is organized as follow.
In Section 2, some basic concepts and useful fundamental results about linear relations are introduced.
In Section 3, we first show that the deficiency indices of Hermitian relations are invariant under relatively bounded perturbations
with relative bounds less than 1.
As a consequence, stability of self-adjointness of Hermitian relations under bounded and relatively bounded perturbations is obtained.
\medskip

\section{Preliminaries}

In this section, we shall recall some basic concepts and introduce some fundamental results about linear relations.

By ${\mathbf C}$ and ${\mathbf R}$ denote the sets of the complex numbers and the real numbers, respectively.
Let $X$ be a complex Hilbert space with inner product $\langle\cdot,\cdot\rangle$.
By $LR(X)$ denote the set of all linear relations of $X\times X$.
Let $T\in LR(X)$. $T$ is said to be a closed relation if $\overline{T}=T$,
where $\overline{T}$ is the closure of ${T}$.
By $CLR(X)$ denote the set of all closed relations of $X\times X$.

Let $S, T\in LR(X)$ and $\alpha \in {\mathbf C}$.
\vspace{-0.2cm}
$$\begin{array}{rrll}
D(T):&=&\{x\in X:\, (x,f)\in T \;{\rm for\; some}\;f\in X\},\\[0.4ex]
R(T):&=&\{f\in X:\, (x,f)\in T \;{\rm for\; some}\;x\in X\},\\[0.4ex]
T(x):&=&\{f\in X:\, (x,f)\in T\},\\[0.4ex]
\alpha T:&=&\{(x,\alpha f):\, (x,f)\in T\},\\[0.6ex]
T+S:&=&\{(x,f+g):\,(x,f)\in T, (x,g)\in S\}.\\[0.6ex]
\end{array}\vspace{-0.2cm}$$
The adjoint of $T$ is defined by\vspace{-0.2cm}
$$T^*=\{(y,g)\in X^2:\,\langle g,x\rangle=\langle y,f\rangle\;{\rm
for\;all}\; (x,f)\in T\}.\vspace{-0.2cm}$$
$T\in LR(X)$ is called a Hermitian relation if $T\subset T^*$,
and it is called a self-adjoint relation if $T=T^*$.
\medskip

\noindent{\bf Definition 2.1 {\rm [11, Definition 2.3]}.}
Let $T\in LR(X)$. The subspace $R(T-\lambda I)^\perp$ is called the deficiency space of $T$ and $\lambda$,
and the number $d_\lambda(T):=\dim(R(T-\lambda I))^\perp$ is called
the deficiency index of $T$ and $\lambda$.
\medskip

Let $T\in LR(X)$ be a Hermitian relation. By [11, Theorem 2.3], $d_\lambda(T)$ is constant
in the upper and lower half-planes; that is, $d_\lambda(T)=d_+(T)$ for all $\lambda\in \mathbf{C}$ with ${\rm Im}\lambda>0$
and $d_\lambda(T)=d_-(T)$ for all $\lambda\in \mathbf{C}$ with ${\rm Im}\lambda<0$, where $d_{\pm}(T):=d_{\pm i}(T)$.
The pair $(d_+(T), d_-(T))$ is called the deficiency indices of $T$, and $d_+(T)$ and $d_-(T)$ are
called the positive and negative deficiency indices of $T$, respectively.
\medskip

In the following, we shall recall concepts of the norm of a subspace and relatively boundedness of two subspaces, and their fundamental properties.

Let $E$ be a closed subspace of $X$. Define the following quotient space [6]:
$$X/E:=\{[x]:\,x\in X\},\quad [x]=\{x\}+E.$$
We define an inner product on the quotient space $X/E$ by
$$\langle [x],[y] \rangle:=\langle x^\bot, y^\bot\rangle,\quad  [x],[y]\in X/E,                                \eqno(2.1)$$
where $x=x_0+x^\bot$, $y=y_0+y^\bot$ with $x_0,\,y_0\in E$ and $x^\bot,\,y^\bot\in E^\bot.$
It can been easily verified that $X/E$ with this inner product is a Hilbert space.
The norm induced by this inner product is the same as that of $X/E$ induced by the norm of $X$.

Now, define the following natural quotient map:
$$Q_E^X:\:X \rightarrow X/E,\:x \mapsto [x].$$
Let $T\in LR(X)$.
By $Q_T$ denote $Q_{\overline{T(0)}}^X$ for briefness without confusion. Define
$$\tilde{T}_s=G(Q_T)T.                                                                                     \eqno(2.2)$$
Then $\tilde{T}_s$ is a linear operator with domain $D(T)$ [3, Proposition II.1.2].
The norm of $T$ at $x\in D(T)$ and the norm of $T$ are defined by, respectively (see [3, II.1]),
$$\begin{array}{rrll}
&&\|T(x)\|:=\|\tilde{T}_s(x)\|,\\
&&\|T\|:=\|\tilde{T}_s\|=\sup\{\|\tilde{T}_s(x)\|:\,x\in D(T)\, with\, \|x\|\leq 1\}.
\end{array}                                                                                                 \eqno(2.3)$$

\noindent\textit{{\bf Lemma 2.1 {\rm ([3, Propositions II.1.4-II.1.7])}.}
Let $S, T\in LR(X)$. Then
\begin{itemize}\vspace{-0.2cm}
\item[{\rm (1)}] $\|Tx\|=d(T(x),0)=d(T(x),T(0))=d(y,\overline{T(0)})=d(y,T(0))$ for every $x\in D(T)$ and $y\in T(x)$;\vspace{-0.2cm}
\item[{\rm (2)}] $\|(\alpha T)(x)\|=|\alpha|\|T(x)\|$,\quad $\|\alpha T\|=|\alpha|\|T\|$ for every $x\in D(T)$ and $\alpha\in \mathbf{C}$;\vspace{-0.2cm}
\item[{\rm (3)}] $\|S(x)+T(x)\| \leq \|S(x)\|+\|T(x)\|$ for $x\in D(T)\cap D(S)$, \quad $\|S+T\| \leq \|S\|+\|T\|$.\vspace{-0.2cm}
\end{itemize}}
\medskip

\noindent{\bf Definition 2.2 {\rm [3, Definition VII.2.1]}.} Let $S, T\in LR(X)$.
\begin{itemize}\vspace{-0.2cm}
\item [{\rm (1)}] $S$ is said to be $T$-bounded if $D(T)\subset D(S)$ and
there exists a constant $c\geq 0$ such that
$$\|S(x)\|\leq c(\|x\|+\|T(x)\|),\quad x\in D(T).                            \vspace{-0.2cm}$$
\item [{\rm (2)}] If $S$ is $T$-bounded, then the infimum of all numbers $b\geq 0$
for which a constant $a\geq 0$ exists such that
$$\|S(x)\|\leq a \|x\|+b\|T(x)\|,\quad x\in D(T),                                                                \eqno(2.4)$$
is called the $T$-bound of $S$.\vspace{-0.2cm}
\end{itemize}

\noindent\textit{\bf Remark 2.1.}
Condition (2.4) is equivalent to the following condition:
$$\|S(x)\|^2\leq a'^2\|x\|^2+b'^2\|T(x)\|^2,\quad x\in D(T),                                                    \eqno(2.5)$$
where the constants $a',b'\geq 0$. It can be easily deduced that (2.5) implies (2.4) with $a=a', b=b'$, whereas (2.4)
implies (2.5) with $a'^2=(1+\varepsilon^{-1})a^2$ and $b'^2=(1+\varepsilon)b^2$ with an arbitrary $\varepsilon>0$.
Consequently, the $T$-bound of $S$ may as well be defined as the infimum of the possible values of $b'$.
\medskip

\noindent\textit{{\bf Lemma 2.2 {\rm ([13, Lemma 2.7])}.} Let $T\in CLR(X)$ and $\lambda\in \mathbf{C}$.
If there exists $c>0$ such that $\|f-\lambda x\|\geq c\|x\|$ for any $(x, f)\in T$,
then $R(T-\lambda I)$ is closed.}\medskip

\noindent\textit{{\bf Lemma 2.3 {\rm ([15, Propositions 2.1, 3.1, 3.3, Theorem 6.3])}.} Let $S, T\in LR(X)$.
\begin{itemize}\vspace{-0.2cm}
\item [{\rm (1)}] $T=T-S+S$ if and only if $D(T)\subset D(S)$ and $S(0)\subset T(0)$.\vspace{-0.2cm}
\item [{\rm (2)}] If $S$ is $T$-bounded with $T$-bound less than $1$,
then $T+S$ is closed if and only if $T$ is closed.\vspace{-0.2cm}
\item [{\rm (3)}] If $T$ is Hermitian, then $D(T)\subset {T(0)}^\perp$ and
$$\langle \tilde{T}_s(x_2),[x_1] \rangle=\langle [x_2], \tilde{T}_s(x_1) \rangle,\,\,x_1, x_2\in D(T).$$\vspace{-0.2cm}
\end{itemize}}
\vspace{-0.6cm}

\noindent\textit{{\bf Lemma 2.4 {\rm ([5, Lemma 2.5])}.} Let $T\in LR(X)$ be Hermitian.
If there is $\lambda\in \mathbf{C} \backslash \mathbf{R}$ such that
$R(T-\lambda I)=R(T-\bar{\lambda} I)=X$, then $T$ is a self-adjoint relation.}\medskip

\noindent\textit{{\bf Lemma 2.5 {\rm ([16, Lemma 5.8])}.} Let $T\in LR(X)$ be self-adjoint.
If $S\in LR(X)$ is Hermitian and $D(T)\subset D(S)$, then $S(0)\subset T(0)$.}\medskip

\section{ Main results}

In this section, we shall first study the stability of deficiency indices of Hermitian relations under perturbations.
Then, we shall use these results to study the stability of self-adjointness of Hermitian relations. \medskip

We shall first prove some useful lemmas:\medskip

\noindent\textit{{\bf Lemma 3.1.}
Let $T\in CLR(X)$ and satisfy that there exists $c>0$ such that
$$\|T(x)\|\geq c\|x\|, \,\,\forall x\in D(T).                                                                        \eqno(3.1)$$
Let $S\in LR(X)$ with $D(T)\subset D(S)$ and $S(0)\subset T(0)$, and satisfy
$$\|S(x)\|\leq b\|T(x)\|, \,\,\forall x\in D(T),                                                                     \eqno(3.2)$$
for some constant $0\leq b<1$.
Then $T+S$ is closed and satisfies (3.1) with $c_1=(1-b)c$ instead of $c$.
Moreover,
$$\dim R(T+S)^\perp=\dim R(T)^\perp.                                                                                 \eqno(3.3)$$
}
\vspace{-5mm}

\noindent{\bf Proof.}
Note that $T$ is closed. It follows from (3.2) and (2) of Lemma 2.3 that $T+S$ is closed.
By Lemmas 2.1 and 2.3, (3.1), and (3.2), we have that for any $(x,f)\in T$ and $(x,g)\in S$,
$$\|f+g\|\geq \|(T+S)(x)\|\geq \|Tx\|-\|Sx\|\geq (1-b)\|Tx\|\geq (1-b)c\|x\|.                                      \eqno(3.4)$$
This implies that $T+S$ satisfies (3.1) with $c_1=(1-b)c$ instead of $c$.

In addition, by (3.1), (3.4), and the closedness of $T$ and $T+S$, $R(T+S)$ and $R(T)$ are closed by Lemma 2.2.

Now we show that (3.3) holds.
Suppose that $\dim R(T+S)^\perp<\dim R(T)^\perp$.
Then there exists $h\in R(T)^\perp\bigcap R(T+S)$ with $h\neq 0$.
Thus, there exists $x\in D(T)$ such that $(x,h)\in T+S$.
And there exist $f\in T(x)$ and $g\in S(x)$ such that $h=f+g$.
In addition, since $T$ is closed, it can be easily verified that $T(0)$ is closed.
Hence, $X/\overline{T(0)}=X/{T(0)}$.
For clarity, for every $z\in X$, by $[z]_T$ denote the element of $X/{\overline{T(0)}}$.
By the assumption that $S(0)\subset T(0)$, it follows from Lemma 2.1 that
$$\|[z]_T\|=d(z,T(0))\leq d(z,S(0))=\|[z]_S\|,\quad z\in X.                                                   \eqno(3.5)$$
Since $T(0)$ is closed, there exist $f_{0}\in T(0)$ and $f^\perp \in T(0)^\perp$ such that $f=f_{0}+f^\perp$.
Further, by noting that $h\in R(T)^\perp \subset T(0)^\perp$ and $f\in R(T)$, it follows that
$$\langle [f]_T,[h]_T \rangle=\langle f^\perp,h \rangle=\langle f,h \rangle=0.$$
It follows that
$$\|T(x)\|^2=\|\tilde{T}_s(x)\|^2=\langle [f]_T,[f]_T \rangle=-\langle [h]_T-[f]_T,[f]_T \rangle
=-\langle [g]_T,[f]_T\rangle,$$
which, together with (3.2) and (3.5), yields that
$$\|T(x)\|^2\leq \|[g]_T\|\|[f]_T\|\leq \|[g]_S\|\|[f]_T\|=\|S(x)\|\|T(x)\|\leq b\|T(x)\|^2.                 \eqno(3.6)$$
We claim that $\|T(x)\|\neq 0$.
If $\|T(x)\|=0$, then $\|S(x)\|=0$ by (3.2).
Hence, $[f]_T=[g]_S=[0]_T$.
This implies that $f\in T(0)$ and $g\in \overline{S(0)}\subset T(0)$.
Consequently, $h=f+g\in T(0)$.
Further with $h\in R(T)^\perp$, we have $h=0$. This is a contradiction with the assumption that $h\neq 0$.
Therefore, $\|T(x)\|\neq 0$.
In view of $0\leq b<1$, it follows from (3.6) that
$$\|T(x)\|^2<\|T(x)\|^2.$$
This is a contradiction.
Hence, $\dim R(T+S)^\perp \geq \dim R(T)^\perp$.

On the other hand, if $\dim R(T+S)^\perp>\dim R(T)^\perp$,
we can similarly find $f\in R(T+S)^\perp\bigcap R(T)$ with $f\neq 0$,
and there exists $x\in D(T)$ such that $(x,f)\in T$.
Set $(x,g)\in S$. Then $(x,f+g)\in T+S$, and consequently $f+g\in R(T+S)$.
By the assumption that $S(0)\subset T(0)$, we have that $(T+S)(0)=T(0)$. Then
$$\|[z]_{T+S}\|= \|[z]_{T}\| \leq \|[z]_S\|,\quad z\in X.                                                     \eqno(3.7)$$
In addition, in view of that $(T+S)(0)$ is closed, there exist $(f+g)_0\in (T+S)(0)$ and $(f+g)^\perp \in ((T+S)(0))^\perp$
such that $f+g=(f+g)_0+(f+g)^\perp$.
Noting that $(T+S)(0)\subset R(T+S)$, we get that $R(T+S)^\perp \subset ((T+S)(0))^\perp$.
So, $f\in R(T+S)^\perp\subset((T+S)(0))^\perp$.
It follows that
$$\langle  [f]_{T+S},[f+g]_{T+S} \rangle=\langle f, (f+g)^\perp \rangle=\langle f, f+g\rangle=0,$$
which, together with (3.7), yields that
$$\|T(x)\|^2=\|[f]_T\|^2=\|[f]_{T+S}\|^2=\langle [f]_{T+S},[f+g]_{T+S}-[g]_{T+S}\rangle
=-\langle [f]_{T+S},[g]_{T+S} \rangle,$$
which implies that
$$\|T(x)\|^2\leq \|[f]_{T+S}\|\|[g]_{T+S}\|\leq \|[f]_T\|\|[g]_S\|=\|T(x)\|\|S(x)\|\leq b\|T(x)\|^2,$$
in which (3.2) has been used.
This is a contradiction with that $0\leq b<1$ and $\|T(x)\|\neq 0$ by $f\neq 0$.
Therefore, $\dim R(T+S)^\perp\leq\dim R(T)^\perp$.
And consequently, $\dim R(T+S)^\perp=\dim R(T)^\perp$.
The whole proof is complete.
\medskip

\noindent\textit{{\bf Corollary 3.1.} Let $T\in CLR(X)$ and satisfy (3.1) for some constant $c>0$.
Let $S\in LR(X)$ with $D(T)\subset D(S)$ and $S(0)\subset T(0)$, and satisfy
$$\|S(x)\|\leq a\|x\|,\: x\in D(T),                                                                      \eqno(3.8)$$
where $0\leq a<c$. Then all the conclusions of Lemma 3.1 hold.}
\medskip

\noindent{\bf Proof.}
It follows from (3.1) and (3.8) that
$$\|S(x)\| \leq a \|x\|\leq \frac{a}{c}\|T(x)\|,\: x\in D(T).$$
Hence, (3.2) is satisfied with $b=\frac{a}{c}<1$.
Therefore, the assertion holds by Lemma 3.1.
The proof is complete.
\medskip

\noindent\textit{{\bf Lemma 3.2.} Let $T\in CLR(X)$ be Hermitian.
Then for any $x\in D(T)$ and any $z\in \mathbf{C}$ with $z=a+ib$ and $a,\,b\in\mathbf{R}$,
$$\|(T-zI)(x)\|^2=\|(T-aI)(x)\|^2+b^2\|x\|^2.                                                            \eqno(3.9)$$}
\vspace{-6mm}

\noindent{\bf Proof.}
Fix any $x\in D(T)$ and any $z=a+ib$ with $a,\,b\in\mathbf{R}$. We have that
$$\begin{array}{rrll}
&&\|(T-zI)(x)\|^2=\|Q_T(T-zI)(x)\|^2\\
&=&\|Q_TT(x)-aQ_T(x)-ibQ_T(x)\|^2\\
&=&\langle \tilde{T}_s(x)-a[x]-ib[x],\tilde{T}_s(x)-a[x]-ib[x]\rangle,
\end{array}                                                                                             \eqno(3.10)$$
where $[x]\in X/\overline{T(0)}$.
Since $T$ is Hermitian, it follows from (3) of Lemma 2.3 that $x\in T(0)^\perp$ and
$$\langle \tilde{T}_s(x),[x] \rangle=\langle [x],\tilde{T}_s(x) \rangle,$$
which implies that
$$\langle \tilde{T}_s(x)-a[x],-ib[x] \rangle=\langle ib[x],\tilde{T}_s(x)-a[x] \rangle.$$
Inserting it into (3.10), we get that
$$\|(T-zI)(x)\|^2=\|\tilde{T}_s(x)-a[x]\|^2+b^2\|[x]\|^2=\|(T-aI)(x)\|^2+b^2\|[x]\|^2.                   \eqno(3.11)$$
In addition, $\|[x]\|=\|x\|$ by noting that $x\in T(0)^\perp$.
Therefore, it follows from (3.11) that (3.9) holds.
This completes the proof.
\medskip

\noindent\textit{\bf Remark 3.1.}
Lemma 3.2 extends the result given in [6, p. 270] for closed symmetric operators to closed Hermitian relations.
\medskip

\noindent\textit{{\bf Lemma 3.3.}
Let $S, T\in LR(X)$ be Hermitian relations with $D(T)\subset D(S)$ and $S(0)\subset T(0)$.
Suppose that $T$ is closed and $S$ is $T$-bounded with $T$-bound less than 1.
Then $T+S$ is closed and $d_\pm(T+S)=d_\pm(T)$.}
\medskip

\noindent{\bf Proof.}
By the assumption that $S$ is $T$-bounded with $T$-bound less than 1, there exist $a>0$ and $0<b<1$ such that
$$\|S(x)\|\leq a\|x\|+b\|T(x)\|, \,\,x\in D(T).$$
Since $0<b<1$, there exists $\varepsilon>0$ such that $0<(1+\varepsilon)b^2<1$.
By Remark 2.1 we have that
$$\|S(x)\|^2\leq (1+\varepsilon^{-1})a^2\|x\|^2+(1+\varepsilon)b^2\|T(x)\|^2, \,\,x\in D(T).$$
Let $\gamma=\frac{a}{b\sqrt{\varepsilon}}$. Then $\gamma>0$ and
$$\|S(x)\|^2\leq (1+\varepsilon)b^2(\gamma^2\|x\|^2+\|T(x)\|^2),\,\,x\in D(T).                            \eqno(3.12) $$
In addition, since $T$ is Hermitian, by Lemma 3.2 we get that
$$\|(T\pm i\gamma I)(x)\|^2=\|T(x)\|^2+\gamma^2\|x\|^2,\,\,x\in D(T).$$
This, together with (3.12), yields that
$$\|S(x)\|\leq (1+\varepsilon)^{1/2}b\|(T\pm i\gamma I)(x)\|,\,\,x\in D(T).$$
Since $0<(1+\varepsilon)^{1/2}b<1$, applying Lemma 3.1 to $T\pm i\gamma I$ and $S$ we get that
$T+S$ is closed and $d_\pm(T+S)=d_\pm(T)$.
This completes the proof.
\medskip

Now, we give the main result of the present paper.\medskip

\noindent\textit{{\bf Theorem 3.1.}
Let $S, T\in LR(X)$ be Hermitian relations with $D(T)\subset D(S)$ and $S(0)\subset T(0)$.
If $S$ is $T$-bounded with $T$-bound less than 1, then $T+S$ is self-adjoint if and only if $T$ is self-adjoint.}
\medskip

\noindent{\bf Proof.}
Since $S$ is $T$-bounded with $T$-bound less than 1, it follows from (2) of Lemma 2.3 that
$T+S$ is closed if and only if $T$ is closed.
In this case, by Lemma 3.3 we get that $d_\pm(T+S)=d_\pm(T)$.
This, together with Lemma 2.4, yields that $T+S$ is self-adjoint if and only if $T$ is self-adjoint.
The proof is complete.
\medskip

The following result is a direct consequence of Lemma 2.5 and Theorem 3.1.
This result is the same as that of [16, Theorem 5.2].
\medskip

\noindent\textit{{\bf Corollary 3.2.}
Let $S\in LR(X)$ be Hermitian and $T\in LR(X)$ be self-adjoint with $D(T)\subset D(S)$.
If $S$ is $T$-bounded with $T$-bound less than 1, then $T+S$ is self-adjoint.}
\medskip

\noindent\textit{\bf Remark 3.2.}
By the definition of relative boundedness for subspaces,
we shall remark that the results about stability of deficiency indices and self-adjointness obtained in the present paper still hold under bounded perturbations.
\medskip

%
%

%
\section*{Funding}
\setcounter{equation}{0}
This work is supported by the NNSF of China (Grant 11571202) and the NSF of Jiangsu Province, China (Grant BK20170298).
%
%
\section*{Competing interests}
\setcounter{equation}{0}
The author declares that she has no competing interests regarding the publication of this paper.

%

\end{document}